\documentclass{birkjour}
\numberwithin{equation}{section}
 \newtheorem{theorem}{Theorem}[section]
\newtheorem{corollary}[theorem]{Corollary}
\newtheorem{lemma}[theorem]{Lemma}
\newtheorem{remark}[theorem]{Remark}
\newtheorem{definition}[theorem]{Definition}
\newtheorem{proposition}[theorem]{Proposition}

\begin{document}

%
%
%
%
%
%
%
%
%

\title[The first eigenvalue of the $p$-Laplacian on time dependent...]
{The first eigenvalue of the $p$-Laplacian on time dependent Riemannian metrics}

\author[Abimbola Abolarinwa]{Abimbola Abolarinwa}

\address{%
Department of Mathematics and Statistics,\\
Osun State College of Technology,\\
P.M.B. 1011, Esa-Oke,\\
 Nigeria.}
 
\email{A.Abolarinwa1@gmail.com}



\author[Jing Mao]{Jing Mao}
\address{%
Department of Mathematics,\\
Harbin Institute of Technology(Weihai),\\
Weihai, 264209,\\
 China.}
\email{jiner120@163.com}

\subjclass{Primary 53C21, 53C44;  Secondary 58C40}

\keywords{Geometric flows, $p$-Laplacian, eigenvalues, monotonicity}

\date{January 1, 2004}

\begin{abstract}
 Let $(M,g)$ be an $n$-dimensional compact Riemannian manifold ($n>1$)  whose metric $g(t)$ evolves by the generalized abstract geometric flow. This paper discusses the evolution, monotonicity and differentiability for the first eigenvalue of the $p$-Laplacian on $(M,g(t))$ with respect to time evolution. We prove that the first nonzero $p$-eigenvalue is monotone nondecreasing along the flow under certain geometric condition and that the first eigenvalue is differentiable 
almost everywhere. When $p=2$, we recover the corresponding results for the usual Laplace-Beltrami operator. Our results provide a unified approach to the study of $p$-eigenvalue under various geometric flows.\end{abstract}

\maketitle


 \section{Introduction} 
 Let $g(t)$ be a one parameter family of Riemannian metrics for $ t \in [0, T]$.
  We say that the time-dependent $g_{ij}(x,t)$ is a  generalized geometric flow if it is evolving by the following 
  equation
\begin{equation}\label{eq11}
\frac{\partial }{ \partial t} g_{ij} (x,t)  =  - 2 h_{ij}(x,t),  \ \ \ \ (x, t) \in M \times [0,T]
\end{equation}
with $g_{ij}( x,0) = g_0(x)$, where $x \in M$, $0<T<T_\epsilon$ is the maximal time of existence, i.e., $T_\epsilon$ is the first time where the flow blows-up and $h_{ij}$ is a general time-dependent symmetric $(0, 2)$-tensor. Here $h_{ij}$ is smooth in both variables $t$ and $x$. This is obvious since $g_{ij}$ is smooth in both variables also.
The scaling factor $2$ in (\ref{eq11}) is insignificant while the negative sign may be important
 in some specific applications for the purpose of keeping the flow either forward or backward in time. Two popular examples of geometric flows are the Ricci flow where $h_{ij} = R_{ij}$ (the Ricci curvature tensor) and the mean curvature flow, in which case $h_{ij} = H \Pi_{ij}$ (here $H$ is the mean curvature and $\Pi_{ij}$ is the second fundamental form on $M$).
 
 One may impose boundedness condition on the tensor $h_{ij}$. In fact, such  boundedness and sign assumptions on $h_{ij}$ are preserved  as long as the flow exists, so it follows that the metrics are uniformly equivalent. Precisely, if $ -K_1 g \leq h \leq K_2g$, where $g(t), t \in [0, T]$ is the flow, then 
\begin{equation}\label{eq2.3}
e^{-K_1 T}g(0) \leq g(t) \leq e^{K_2 T}g(0).
\end{equation}
To see the above bounds (\ref{eq2.3}) we consider the evolution of a vector form $|X|_g = g(X, X), X \in T_xM$. By the equation (\ref{eq11}) and by the boundedness of the tensor $h$ we have 
$ | \partial_t g(X, X) | \leq K_2 g(X, X),$
which implies (by integrating from $t_1\ to\ t_2$)
$$ \Big| \log \frac{g(t_2)(X, X)}{g(t_1)(X, X)}\Big| \leq K_2 t \Big|_{t_1}^{t_2}.$$
Taking the exponential of this estimate with $t_1 =0$ and $t_2 =T$ yields $|g(t)| \leq e^{k_2 T}g(0)$ from which the uniform boundedness of the metric follows. Thus, if there holds boundedness assumption $ -K_1 g \leq h \leq K_2g$, the metric $g(t)$ are uniformly bounded below and above for all time $0 \leq t \leq T$ under the geometric flow (\ref{eq11}). Then, it does not matter what metric we use in the argument that follows.

 In this paper, we 
 discuss the evolution, monotonicity and differentiability of the first nonzero eigenvalue of $p$-Laplace operator on the $n$-dimensional compact Riemannian manifold whose metric evolves by (\ref{eq11}). The $p$-Laplace operator is defined as 
 $$ \Delta_{p, g} f(x) := \ div  (|\nabla f|^{p-2} \nabla f)(x)$$
for $ p \in [1, \infty)$, where $div$ is the divergence operator, the adjoint of gradient ($grad$) for the $L^2$-norm induced by $g$ on the space of differential forms. When $p=2$,
$$ \Delta_{2, g} f(x) = div \circ grad \ f(x) $$
is the usual Laplace-Beltrami operator. The eigenvalues and the corresponding eigenfunctions of $ \Delta_{p, g}$ satisfy the following nonlinear eigenvalue problem
\begin{equation}\label{eq12}
 \Delta_{p, g} f = - \lambda |f|^{p-2} f , \hspace{1cm} f \neq 0.
\end{equation}

It is easily verifiable that the principal symbol of (\ref{eq12}) is nonnegative everywhere and  strictly positive at the neighbourhood of the point where $\nabla f \neq 0$.
We know that (\ref{eq12}) has weak solutions with only partial regularity of class $C^{1, \alpha}, (0< \alpha < 1)$ in general. Interested readers can find the classical papers by L. Evans \cite{[Ev81]} and P. Tolksdorff \cite{[Tok84]}. Notice that the least eigenvalue of $\Delta_{p,g}$ on compact manifold without boundary is zero with the corresponding eigenfunction being a constant. Hence, we refer to the infimum of the positive eigenvalues as the first nonzero eigenvalue or simply the first eigenvalue. The first eigenvalue of $ \Delta_{p, g}$ is characterised by the min-max principle
\begin{equation}\label{eq13}
  \lambda_{p,1} = \inf_{0\neq f \in W^{1, p}(M) }\Bigg\{ \frac{\int_M | \nabla f |_g^p \ d \mu_g}{\int_M | f|^p_g \ d\mu_g} \ \  \Big| \ \ f \neq 0, \ \ f \in W^{1, p}(M) \Bigg\},
 \end{equation}
satisfying the following constraint $\int_M |f|_g^{p-2} f d \mu_g = 0$, where $d \mu_g$ is the volume measure on $(M,g)$. Obviously, the infimum does not change when  one replaces $ W^{1, p}(M)$ by $C^\infty(M)$.
The corresponding eigenfunction is the energy minimizer of Rayleigh quotient (\ref{eq13}) and satisfies the following Euler-Lagrange equation
\begin{equation}
\int_M [ |\nabla f |^{p-2} \langle \nabla f, \nabla \phi \rangle - \lambda |f|^{p-2}  \langle  f,  \phi \rangle ]d\mu_g= 0
 \end{equation}
for $ \phi \in C^\infty_0(M)$ in the sense of distribution.
The problem of finding $ \lambda_{p,1}$ is related to the problem of finding the best constant $C(M)$ in the $L^p$-Sobolev inequality 
$$ \| f \|_{L^{p^*}} \leq C(M) \| \nabla f \|_{L^p},\ \ \ p^* = np/(n-p),$$
which is obtained by  continuous embedding of $W^{1,p}(M) \hookrightarrow L^p(M)$ under the sobolev norm
$$\|f\|_{W^{1,p}} = \Big(\int_M |f|^p d\mu + \int_M | \nabla f|^p d\mu \Big)^{\frac{1}{p}}.$$

It is well-known that $p$-Laplacian has discrete eigenvalues but still remains unknown whether it only has discrete eigenvalues for bounded connected domains. Another well-known results tell us that the first nonzero eigenvalue is simple and isolated \cite{[Linq1],[Linq2],[MaRos]}. Here the simplicity shows that any nontrivial eigenfunction corresponding to $\lambda_{p,1}$ does not change sign and that any two first eigenfunctions are constant multiple of each other. 

In contrast to the spectrum of  Laplace-Beltrami operator (the case $p=2$), the $p$-Laplacian is nonlinear in general. For instance, the case $p=1$ gives $div(\nabla u /|\nabla u|)$, the negative of mean curvature operator.  Aside being of geometric interests, the $p$-Laplacian appears naturally in the study of non-newtonian fluids, nonlinear elasticity, heat radiation, porous media flow,  rheology, petroleum extraction, Brownian motions, torsional creep problem among others. For details of physical applications of $p$-Laplacian, the reader is referred to \cite{[DiSe]}.

 Moreover, it is not known if $\lambda_{p,1}$ or its corresponding eigenfunction is $C^1$-differentiable (or even locally Lipschitz) along any geometric flow of the form (\ref{eq11}). However, it has been pointed out that the differentiability for the case $p=2$ is a consequence of eigenvalue perturbation theory, see for instance \cite{[Kat84],[KL06]}. For this reason, the method of L. Ma  \cite{[Ma06]}, which assumes differentiability of eigenvalues and eigenfunctions under the Ricci flow, can only be applied to find the monotonicity of first eigenvalue for the case $p=2$ along (\ref{eq11}). Now to avoid the differentiability assumption on the first eigenvalue and the corresponding eigenfunction in the case $p \neq 2$, we shall apply  techniques of Cao  \cite{[Ca07]} as used by Wu \cite{[Wu11]} and Wu, Wang and Zheng \cite{[WWZ]} under the Ricci flow to study the evolution and monotonicity of $\lambda_{p,1}(t) = \lambda_{p,1}(t, f(t))$, where $\lambda_{p,1}(t, f(t))$ and $f(t)$  are assumed to be smooth. The evolution and the monotonicity formulas for the first eigenvalue (in both cases $p=2$ and $p  \neq 2$) derived here do not depend on the evolution of the eigenfunction. The eigenfunction only needs to satisfy certain normalization condition. 

There are many results on the evolution and monotonicity of eigenvalues of the Laplace operator on evolving manifolds with or without curvature assumptions. One can find \cite{[Ca07],[Ca08],[Li07a],[Li07b]} under the Ricci flow, \cite{[Li10]} under Ricci-Harmonic map flow and \cite{[GPT]} along abstract geometric flow with entropy methods. The study of the properties of eigenvalues of the $p$-Laplacian on evolving manifold is still very young. The main aim of this paper is to investigate if those known properties of $\lambda_{p,1}$ on static metric and for the case $p=2$ on evolving metric can be extended to various geometric flow. We however intend to develop a unified algorithm that can be used for this purpose on time-dependent metrics. Many interesting results concerning the behaviour of $\lambda_{p,1}$ can be found in \cite{[KaFr],[Linq1],[Linq2],[Mao2],[MaRos]} for static metrics and \cite{[Ab15],[GLW16],[Mao1],[Wu11],[WWZ],[Zha12],[Zha13]} for evolving metrics along various geometric flow.

The rest of this paper is planned as follows. In Section \ref{sec2}, we discuss the main results of this paper. Firstly, we highlight some notations about differential geometry including metrics, gradient, divergence and integration by parts, which form the core of analytic tools used in the paper. We also prove some technical lemmas about evolutions of some geometric objects relating to the eigenvalues under the flow (\ref{eq11}). In Section \ref{sec3}, we study the first eigenvalue of Laplace-Beltrami operator ($p=2$) under this geometric flow. Here, we assume the first eigenvalue to be a function of time only and obtain its general evolution and monotonicity formula under certain condition. In Section \ref{sec4}, we use a different approach to derive $p$-eigenvalue's evolution and monotonicity without differentiability assumption on the eigenfunction. In fact, the differentiability of  $p$-eigenvalue is a consequence of  the monotonicity formula. In the last section, we list  some examples of geometric flows where the approaches used in this paper are applicable. In fact, this section reveals that our generalised geometric flow is not a trivial generalisation.

\section{Preliminaries and main results}\label{sec2}
\subsection{Notation}
Throughout, $M$ will be taken to be a closed manifold (i.e., compact without boundary).
Most of our calculations are done in local coordinates, where $ \{ x^i \}$ is fixed in a neighbourhood of every point $ x \in M$. We shall adopt Einstein summation convention with repeated indices summed up.
The Riemannian metric $g(x)$ at any point $x \in M$ is a bilinear symmetric positive definite matrix denoted in local coordinates by 
$$g_{ij}(x) =  g_{ij} dx^i dx^j.$$
The Laplace-Beltrami operator acting on a smooth function $ f $ on $M$ is defined as divergence of  gradient of $f$, written as
$$\Delta_g f := div \ grad \ f =\frac{1}{\sqrt{|g|}} \frac{\partial }{\partial x^i} \Big( \sqrt{|g|}g^{ij} \frac{\partial f}{\partial x^j}  \Big), $$
where $|g| = \det (g_{ij})$ and  $g^{ij} = (g_{ij})^{-1}$ are determinant and the inverse metric matrix respectively. By the above we note that 
$$(grad \ f)^i = (\nabla f)^i = g^{ij} \frac{\partial f}{\partial x^j}  \ \ \mbox{and} \ \ div  X = \frac{1}{\sqrt{|g|}} \frac{\partial }{\partial x^i} ( \sqrt{|g|} X^i ), $$
where $X$ is a smooth vector field.
Also we have the metric norm 
$$| \nabla f|^2_g =  g^{ij} \nabla_i f \nabla_j f = \nabla^j f  \nabla_j f.$$
 The Riemann structure  allows us to define Riemannian volume measure $d \mu_g$ on $M$ by
$$ d \mu_g  = \sqrt{|g_{ij}(x)|}dx^i .$$
By the divergence theorem, we have the following integration by parts formulas:
Let $X$ be a vector field, $X = X^i \partial_i$ and $f \in  C^\infty(M)$ be smooth function, then 
$$   \int_M \langle - div X, f \rangle_g  =  \int_M  \langle X, \nabla f \rangle_g = - \int_M \frac{1}{\sqrt {\det g}} f  \partial_i ( X^i \sqrt{\det g }) \sqrt{\det g} \ d x^i .$$
Also for  functions $f, h \in C^2(M)$
$$ \int_M f \Delta_g h \ d\mu = - \int_M \langle \nabla f, \nabla h \rangle_g d\mu = \int_M \Delta_g f \ h d\mu.$$ 
We write in local coordinates gradient, Hessian and covariant derivative as
$$\nabla f = f_i, \ \ \nabla \nabla f = \nabla_i \nabla_j f = f_{ij} \ \  and \ \ \frac{\partial}{\partial x^i} = \partial_i$$
 respectively.
 Also we write time derivative as $\frac{\partial}{\partial t} f = \partial_t f = f_t$. 
 
\subsection{Evolution equations}
Interestingly, all the geometric quantities associated with the underlying manifold evolve as the Riemannian metric evolves along the geometric flow. This also serves as a motivation considering the behaviours of some other important geometric quantities such as eigenvalues  of the manifold under the  flow. The next two lemmas give us these evolutions.
\begin{lemma}\label{lem211}
Suppose a one-parameter family of smooth metrics $g(t)$ solves the geometric flow (\ref{eq11}), then, we have the following evolutions
\begin{align*}   
 \displaystyle & (1) \hspace{2cm}   \frac{\partial}{\partial t } g^{ij} =  2 g^{ik} g^{jl} h_{kl} = 2 h^{ij}\\
 \displaystyle  &  (2)\hspace{2cm}    \frac{\partial}{\partial t } \Gamma^k_{ij} =  -  g^{kl} \Big( \partial_i h_{jl} +   \partial_j h_{il} -    \partial_l h_{ij} \Big)   \\
  \displaystyle &  (3) \hspace{2cm}   \frac{\partial}{\partial t} | \nabla f|^2 = 2 h_{ij} f_i f_j + 2 f_i f_{t,j} \\
   \displaystyle & (4) \hspace{2cm}    \frac{\partial}{\partial t} ( \Delta f ) =  2 h_{ij}f_{ij} +  2 \langle \ div \ h , \nabla f \rangle -   \langle \nabla \mathcal{H}   , \nabla f \rangle + \Delta f_t \\
    \displaystyle & (5) \hspace{2cm}    \frac{\partial}{\partial t}  d \mu = - \mathcal{H} d\mu. 
\end{align*}
\end{lemma}
Here $\mathcal{H} = g^{ij} h_{ij}$, the metric trace of a symmetric $2$-tensor $h_{ij}$ and $f$ is a smooth function defined on $M$.

\begin{proof}
Recall that both $g_{ij}$ and $\mathcal{S}_{ij}$ are symmetric tensors and $ g^{ij}g_{jl} = \delta^i_l$. 
Note  also that Levi-civita connection is not a tensor, but the time derivative of a connection is $(2, 1)$ tensor.  Then $(2)$ holds as a tensor equation in any coordinates system and at any point. The proofs of $(1)$ and $(2)$  are the same as those of \cite[Lemmas 3.1 and 3.2]{[CK04]}. We  give the computations in local coordinates which lead to the proofs of $(3)$ and $(4)$.  In fact, we have
 \begin{align*}
  \partial_t (| \nabla f|^2) & = \partial_t(g^{ij} \partial_i f \partial_j f) \\
 \displaystyle &=  (\partial_t g^{ij}) \partial_i f \partial_j f + 2 g^{ij}   \partial_i f \partial_j   f_t \\  
\displaystyle & =  2 h^{ij}  \partial_i f \partial_j f  + 2 g^{ij}   \partial_i f \partial_j   f_t \\  
\displaystyle & =  2 h_{ij}  f_i f_j   + 2 f_i f_{t, j},
       \end{align*}
 which is exactly (3). 
 We prove (4) by using (1) and (2) as follows:
\begin{align*} 
\frac{\partial}{\partial t} ( \Delta f) & = \frac{\partial}{\partial t} [ g^{ij} ( \partial_i \partial_j - \Gamma_{ij}^k  \partial_k ) f ] \\
& = \frac{\partial}{\partial t} (g^{ij}) ( \partial_i \partial_j - \Gamma_{ij}^k  \partial_k ) f  + g^{ij} ( \partial_i \partial_j - \Gamma_{ij}^k  \partial_k ) \frac{\partial}{\partial t} f  - g^{ij} \Big( \frac{\partial}{\partial t} \Gamma_{ij}^k \Big) \partial_k f \\
& = 2 h^{ij} \nabla_i \nabla_j f + \Delta f_t + g^{ij} g^{kl} \Big( \nabla_i h_{jl}  + \nabla_j h_{il} - \nabla_l h_{ij}\Big) \nabla_k f \\
&=  \Delta f_t + 2 h_{ij} f_{ij} + 2 g^{kl}  \Big( g^{ij} \nabla_i h_{jl} - \frac{1}{2} g^{ij} \nabla_l h_{ij} \Big) \nabla_k f \\
& = \Delta f_t + 2 h_{ij} f_{ij} +2 \langle \ div \ h , \nabla f \rangle  - \mathcal{H} _i f_j. 
\end{align*}
In local coordinates, the volume form is written as  $ d \mu = \sqrt{|g| } dx^1 \wedge . . .\wedge dx^n$,  then, 
$$\frac{\partial}{\partial t} d\mu = \frac{\partial}{\partial t} \Big(\sqrt{| g| }dx^1 \wedge . . .\wedge dx^n \Big).$$
 By chain rule of differentiation
  $$\frac{\partial}{\partial t} \Big( \sqrt{| g| } \Big) =   \frac{1}{2} \frac{1}{\sqrt{|g| }} \frac{\partial}{\partial t} | g| =  \frac{1}{2} \frac{1}{\sqrt{|g| }} \frac{\partial | g|}{\partial g_{ij} } \frac{\partial g_{ij}}{\partial t }$$ 
  $$\ \ \ \ \ =  - \sqrt{|g |}g^{ij}h_{ij} = - \mathcal{H}  \sqrt{|g| }.$$
    Therefore
     $$\frac{\partial}{\partial t} d \mu = - \mathcal{H} \ d\mu.$$
This completes the proof of the lemma.
\end{proof}

\begin{lemma}\label{lem212}
Suppose a one-parameter family of smooth metrics $g(t)$ solves the geometric flow (\ref{eq11}), then, we have the following evolutions
\begin{align*}   
  \displaystyle &  (1) \hspace{2cm}   \frac{\partial}{\partial t} | \nabla f|^p = p | \nabla f|^{p-2} \Big\{ h_{ij} \nabla_i f  \nabla_j  f + g^{ij} \nabla_i f \nabla_j  f_t   \Big\}\\
   \displaystyle &  (2) \hspace{2cm}   \frac{\partial}{\partial t} | \nabla f|^{p-2} = (p-2) | \nabla f|^{p-4} \Big\{ h_{ij} \nabla_i f  \nabla_j  f + g^{ij} \nabla_i f \nabla_j  f_t   \Big\}\\
   \displaystyle & (3) \hspace{2cm}    \frac{\partial}{\partial t} ( \Delta_{p, g} f ) =   2 h^{ij} \nabla_i ( Z \nabla_j f ) + g^{ij} \nabla_i ( Z_t \nabla_j f ) + g^{ij} \nabla_i (Z \nabla_j f_t )\\
  \displaystyle & \hspace{5cm}   +  Z \Big\{2 \langle \ div \ h , \nabla f \rangle -   \langle \nabla \mathcal{H}   , \nabla f \rangle + \Delta f_t \Big\},
  \end{align*}
  where $Z:= | \nabla f|^{p-2}$ and $f$ is a smooth function defined on $M$. When $p=2$ we have (2) and (3) of Lemma \ref{lem211}.
\end{lemma}

The proof follows standard computation as in \cite[Lemma 2.2]{[Ab15]}. But we include it here for completeness. 
\begin{proof}
By (3) of Lemma \ref{lem211} we can get
\begin{align*}
\frac{\partial}{\partial t} \Big(| \nabla f|^p \Big) & =  \frac{\partial}{\partial t} \Big( | \nabla f|^2\Big)^{\frac{p}{2}}\\
   \displaystyle &=  \frac{p}{2}  \Big( | \nabla f|^2\Big)^{\frac{p-2}{2}} \frac{\partial}{\partial t} \Big( | \nabla f|^2\Big)\\
  \displaystyle &=  \frac{p}{2}   | \nabla f|^{p-2} \Big\{  2 h_{ij} \nabla_i f  \nabla_j  f + 2 g^{ij} \nabla_i f \nabla_j  f_t  \Big\} \\  
\displaystyle & =  p  | \nabla f|^{p-2} \Big\{   h_{ij} \nabla_i f  \nabla_j  f +  g^{ij} \nabla_i f \nabla_j  f_t  \Big\},
       \end{align*}
which is (1). (2) follows from the same calculation as the above. Let $Z:= | \nabla f|^{p-2}$, then 
 \begin{align*}
\frac{\partial}{\partial t} \Delta{p, g} f  & =  \frac{\partial}{\partial t} \Big(div ( | \nabla f|^{p-2} \nabla f )\Big) = \frac{\partial}{\partial t} \Big(g^{ij} \nabla_i (Z \nabla_j f )\Big) \\
   \displaystyle &=  \frac{\partial}{\partial t} \Big(g^{ij} \nabla_i Z \nabla_j f  + Z g^{ij} \nabla_i \nabla_j f \Big) \\
   \displaystyle &=   \frac{\partial}{\partial t} \Big( g^{ij} \Big) \nabla_i Z \nabla_j f + g^{ij} \nabla_i Z_t \nabla_j f + g^{ij} \nabla_i Z \nabla_j f_t + Z_t \Delta f + Z (\Delta f)_t .
       \end{align*}
       By (1) and (3) of the Lemma, we have
     \begin{align*}
\frac{\partial}{\partial t} \Delta{p, g} f  & = 2 h^{ij} \nabla_i Z \nabla_j f + g^{ij} \nabla_i Z_t \nabla_j f + g^{ij} \nabla_i Z \nabla_j f_t + Z_t \Delta f \\
  \displaystyle &\hspace{1cm} + Z \Big\{ 2 h_{ij} \nabla_i \nabla_j f +  2 \langle \ div \ h , \nabla f \rangle -   \langle \nabla \mathcal{H}   , \nabla f \rangle + \Delta f_t \Big\}\\
 \displaystyle & =  2 h^{ij} \nabla_i Z \nabla_j f  + 2 h_{ij} Z  \nabla_i \nabla_j f + g^{ij} \nabla_i Z_t \nabla_j f  + Z_t \Delta f  + g^{ij} \nabla_i Z \nabla_j f_t\\
  \displaystyle &\hspace{1cm} + Z \Delta f_t + Z \Big\{2 \langle \ div \ h , \nabla f \rangle -   \langle \nabla \mathcal{H}   , \nabla f \rangle \Big\}\\  
   \displaystyle & =  2 h^{ij} \nabla_i ( Z \nabla_j f ) + g^{ij} \nabla_i ( Z_t \nabla_j f ) + g^{ij} \nabla_i (Z \nabla_j f_t ) \\
   \displaystyle &\hspace{1cm} +  Z \Big\{2 \langle \ div \ h , \nabla f \rangle -   \langle \nabla \mathcal{H}   , \nabla f \rangle \Big\}.
 \end{align*}
   \end{proof}

\subsection{Main Results}

Recall that we already mentioned that sign assumption on the tensor $h_{ij}$ is preserved throughout the flow. To prove the monotonicity of $\lambda_{p,1}$, we will need the condition $(h_{ij} - \alpha \mathcal{H} g_{ij})(x,t) \geq 0$ for all $t \in [0, T]$. This condition is informed by the Hamilton's maximum principle for tensors. For clarity we state the principle without proof.
\begin{definition}({\bf Null-eigenvector assumption})
A quantity $Q(q,t) : Sym^2T^*M \times [0,T] \to Sym^2T^*M$ is said to satisfy the null eigenvector assumption if whenever $\omega_{ij}$ is a nonnegative symetric $(0,2)$-tensor at a point $q$ and if $X \in T_qM$ is such that $\omega_{ij} X^j =0$ and then 
$$Q_{ij}(\omega,g)X^iX^j \geq 0$$
for any $t \in [0,T]$.
\end{definition}
The symetric tensor $\omega_{ij}$ is defined to be nonnegative if and only if $\omega_{ij} V^iV^j \geq 0$ for all vectors $V^i$ (i.e., the quadratic form induced by $\omega_{ij}$ is semi-positive definite) and we write $\omega_{ij} \geq 0$.

\begin{theorem}(\cite[Theorem 9.1]{[Ha82]},\cite[Theorem 4.6]{[CK04]})
Let $g(t)$ be a smooth one parameter family of Riemannian metrics satisfying (\ref{eq11}). Let $M_{ij}$ be a symmetric $(0,2)$-tensor satisfying
\begin{align}
\frac{\partial}{\partial t} M_{ij}(x,t) \geq \Delta M_{ij}(x,t) + \langle X, \nabla M_{ij}\rangle + Q(M_{ij},g(t)),
\end{align}
where $X$ is a time-dependent vector field and $Q(M_{ij},g(t))$ is a symmetric $(0,2)$-tensor which is locally Lipschitz in $x$, continuos in $t$ and satisfies the null-eigenvector assumption. 

If $M_{ij}(p,0) \geq 0$ for all $p \in M$, then $M_{ij}(p,t) \geq 0$ for all $p \in M$ and $t \in [0,T_\epsilon)$
\end{theorem}

By the above maximum principle for tensors we can prove the following conclusion.
\begin{proposition}
Let $g(t)$ be a smooth one parameter family of Riemannian metrics satisfying (\ref{eq11}). If 
$$(h_{ij} - \alpha \mathcal{H} g_{ij})(x,0) \geq 0,$$
then
\begin{align}
(h_{ij} - \alpha \mathcal{H} g_{ij})(x,t) \geq 0
\end{align}
for some $\alpha \in [0,\frac{1}{n}]$ and all $t \in [0, T].$
\end{proposition}

The first main results on evolution and monotonicity of $\lambda_{p,1}$ for the case $p=2$ are the following and are proved in Section \ref{sec3}.
\begin{theorem}\label{thm311}
Let $(M, g(t))$ be a closed manifold evolving by the geometric flow (\ref{eq11}) and $\lambda(t)$ be the first nonzero eigenvalue of the Laplacian $\Delta_{g(t)}$ corresponding to the eigenfunction $ \varphi(t, x)$, then $\lambda(t)$ evolves by 
\begin{equation}
   \frac{d}{d t} \lambda(t) \int_M   \varphi^2 \ d\mu =  \lambda(t) \int_M  \mathcal{H} \varphi^2 \ d\mu -   \int_M \mathcal{H} | \nabla \varphi |^2 \ d\mu  + 2  \int_M \langle h, d \varphi \otimes d  \varphi \rangle \ d\mu. 
  \end{equation}
\end{theorem}

By setting $ h_{ij} - \alpha \mathcal{H}g_{ij} \geq 0, \alpha \geq 1/2 $, along the flow we have the following monotonicity formula from the last theorem
\begin{equation}
 \frac{d}{d t} \lambda(t) \geq \   \lambda(t)  \frac{\int_M \mathcal{H} \varphi^2 \ d\mu }{ \int_M   \varphi^2 \ d\mu} + ( 2 \alpha - 1) \frac{\int_M \mathcal{H} | \nabla \varphi |^2 \ d\mu }{\int_M   \varphi^2 \ d\mu}.
\end{equation}
By the definition of $\lambda(t)$ or from (\ref{eq21}) we know that $\lambda(t) \int_M \varphi^2 \ d\mu =  \int_M  | \nabla \varphi |^2 \ d\mu$
and $\lambda(t) > 0$. Suppose further that $\mathcal{H}(x,0) \geq 0$, then,  
\begin{equation}\label{eq25}
   \frac{d}{d t} \lambda(t)  \geq 2 \alpha \lambda(t)  \frac{\int_M \mathcal{H} \varphi^2 \ d\mu }{ \int_M   \varphi^2 \ d\mu}  .
  \end{equation}

The monotonicity of $\lambda(t)$ here depends on the sign of  $\mathcal{H}$. Note that in applications the sign of  $\mathcal{H}$ is usually preserved throughout the evolution. An interesting case is when the manifold is being evolved under the Ricci flow \cite{[Ha82]}, where the nonnegativity of scalar curvature is preserved.
Recall also that $\mathcal{H}$  evolves by
  $$\frac{\partial}{ \partial t } \mathcal{H} =  \beta + 2 |h_{ij}|^2$$
  where  $\beta := g^{ij} \partial_t h_{ij}$, in particular,  under Ricci flow where $h_{ij} = R_{ij}$ and $\mathcal{H} = R$, we have $\beta = \Delta R$. Here  we will assume that 
  \begin{equation}\label{eqB}
  \beta - \Delta \mathcal{H} \geq 0.
  \end{equation}
 This is motivated by an error term appearing in a result of M\"uller \cite[Lemma 1.6]{[Mu10]}. For our case the error term reads; for any time-dependent vector field $X$ on $M$
  \begin{equation}\label{eq224}
  \mathcal{D}(X) :=  2(R_{ij} -h_{ij})(X, X) + 2\langle \nabla \mathcal{H} - 2 div \ h , X\rangle + \partial_t \mathcal{H} - \Delta \mathcal{H} - 2|h_{ij}|^2,
  \end{equation}
 where $R_{ij}$ is the Ricci curvature tensor of $M$. Clearly the last three terms in (\ref{eq224}) above is the same as the quantity $ \beta - \Delta \mathcal{H}$. It does make sense to assume (\ref{eqB}) holds whenever $ \mathcal{D}(X)$ is nonnegative. An application of this is that we are on a steady or shrinking  soliton (self-similar solution to the geometric flow) if the equality in (\ref{eqB}) holds.
Writing  $|h_{ij}|^2 \geq \frac{1}{n} \mathcal{H}^2$ and  using the condition that $\beta - \Delta \mathcal{H} \geq 0$, we have a governing differential inequality for the evolution of $\mathcal{H}$ as follows
\begin{equation}\label{eq318}
\frac{\partial}{\partial t} \mathcal{H} \geq \Delta \mathcal{H} + \frac{2}{n} \mathcal{H}^2.
\end{equation}
Suppose $\mathcal{H} \geq \mathcal{H}_{min},$ we can apply the maximum principle by comparing the solution of (\ref{eq318}) with that of the following ordinary differential inequality
 \begin{equation}
    \left \{ \begin{array}{l}
\displaystyle  \frac{d \psi(t) }{d t} = \frac{2}{n} (\psi(t))^2 \\
\displaystyle \  \ \psi(0) = \mathcal{H}_{min}(0),
\end{array} \right.
\end{equation}
solving to
 $$ \psi(t) \ \ \ = \ \ \frac{ \mathcal{H}_{min}(0)}{ 1 - \frac{2}{n} \mathcal{H}_{min}(0)t }.$$  
Therefore  
\begin{equation}\label{Risup}
 \mathcal{H}_{g(t)} \ \ \geq \  \ \psi(t) \  \ =  \  \ \frac{ \mathcal{H}_{min}(0)}{ 1 - \frac{2}{n} \mathcal{H}_{min}(0)t }
  \end{equation}
for all $t \geq 0$ as long as the flow exists.  Hence we write (\ref{eq25}) as follows
\begin{equation}\label{eq29}
   \frac{d}{d t} \lambda(t)  \geq  2 \alpha \psi(t) \lambda(t) .
  \end{equation}
By this we can prove the following
\begin{theorem}
Let $(M, g(t))$ be a closed manifold evolving by the geometric flow (\ref{eq11}). Let $\lambda(t)$ be the first nonzero eigenvalue of the Laplacian $\Delta_{g(t)}$. Suppose
$$ \mathcal{H}_{g(t)} \ \ \geq \  \ \psi(t) \  \ =  \  \ \frac{ \mathcal{H}_{min}(0)}{ 1 - \frac{2}{n} \mathcal{H}_{min}(0)t }$$
Then 
\begin{equation}
   \frac{d}{d t} \Big[ \lambda(t) \exp \Big( - 2 \alpha  \int_0^T   \psi(t) dt \Big) \Big] \geq 0,
  \end{equation}
  where $\alpha \geq \frac{1}{2}$ and 
  \begin{equation}
   \lambda(t) \geq  \lambda(0)  e^{ 2 \alpha  \int_0^T   \psi(t) dt}
  \end{equation}
  for $0 \leq t \leq T$.
\end{theorem}
 
In the next we state our results concerning the general $p$. The proofs are discussed in Section \ref{sec4}.
\begin{theorem}\label{thm411}
Let $(M, g)$ be an $n$-dimensional closed Riemannian manifold evolving by the geometric flow (\ref{eq11}). Let $\lambda_{p,1}(t)$ be the first eigenvalue of the $p$-Laplacian on $M$ corresponding to the eigenfunction $u(t, x)$ at time $t \in [0, T]$.  Then
\begin{equation}\label{eq48}
\left. \begin{array}{l}
\displaystyle\frac{d}{dt} \lambda_{p,1}(t) = \lambda_{p, 1}(t) \int_M \mathcal{H} | u|^p  d\mu -  \int_M \mathcal{H} | \nabla u|^p d\mu \\ 
\displaystyle \hspace{3cm}+ p  \int_M  | \nabla u|^{p-2} h^{ij} \nabla_i u \nabla_j u  d\mu
\end{array} \right. 
\end{equation}
for all time $ t \in[0, T]$.
Moreover, if it holds that 
\begin{equation}\label{eq49}
h_{ij} - \alpha \mathcal{H} g_{ij} \geq 0 , \ \ \ \ \alpha \in [1/p, 1/n),
\end{equation}
then $ \lambda_{p,1}(t)$ is monotonically nondecreasing along the geometric flow and it is differentiable almost everywhere. Precisely,
\begin{equation}\label{eq410}
\frac{d}{dt} \lambda_{p,1}(t) \geq  \lambda_{p, 1}(t) \int_M \mathcal{H} | u|^p  d\mu +(\alpha p -1)  \int_M \mathcal{H} | \nabla u|^p d\mu \geq 0
\end{equation}
provided $\mathcal{H}$ is nonnegative.
\end{theorem}

\begin{corollary}\label{cor4}
With the conditions of Theorem \ref{thm411}.
\begin{equation}\label{eq514}
\lambda_{p, 1}(t_2) \geq \lambda_{p, 1}(t_1) + \int_{t_1}^{t_2} \Theta(g(t), u(t) ) dt,
\end{equation}
where
$$ \Theta(g(t), u(t) ) = \lambda_{p, 1}(t) \int_M \mathcal{H} | u|^p  d\mu +(\alpha p -1)  \int_M \mathcal{H} | \nabla u|^p d\mu.$$
Furthermore, if $\mathcal{H} \geq \mathcal{H}_{min} > 0$ and satisfies the governing inequality (\ref{eq318}), it then holds for all time $t_1 < t_2$ that 
\begin{equation}\label{eq414}
\lambda_{p, 1}(t_2) \geq \lambda_{p, 1}(t_1) \exp\Big\{ \alpha p  \int_{t_1}^{t_2} \psi(t) dt \Big\}.
\end{equation}
\end{corollary}

Finally, we show that the following quantity
\begin{equation}
\lambda_{p,1}(t) \cdot \Big( \psi^{-1}_0 - \frac{2}{n}t \Big)^{\frac{\alpha np}{2}}
\end{equation}
is nondecreasing (see Theorem \ref{thm44})  and that $\lambda_{p,1}(t) $ is differentiable almost everywhere along the geometric flow.

\section{Eigenvalues of the Laplace-Beltrami operator}\label{sec3}
The eigenvalue problem involving Laplace-Beltrami operator on a closed manifold consists in finding all possible real $\lambda$ such that there exists non-trivial functions $u$ satisfyng 
\begin{align}
\Delta_g u = - \lambda u.
\end{align} 
It is well-known that $- \Delta_g$ has a discrete spectrum on a closed Riemannian manifold. This set consists of an infinite sequence
$$0=\lambda_0<\lambda_1 \leq \lambda_2 \leq ...\lambda_k \leq ... \to \infty \ \ \ \mbox{as} \ k \to \infty$$
and can be found that only constant functions correspond to $\lambda_0=0$. The eigenfunctions are $L^2(M,g)$ orthonormal basis $\{u_0, u_1, u_2,...,u_k,...\}$ of real $C^\infty(M,g)$ function such that 
$$\Delta_g u_j = - \lambda_j u_j, \ \ \  j=1,2,...$$
while the eigenvalues are $L^2(M,g)$ orthogonal. By the min-max principle, the first non-zero eigenvalue $\lambda_1(M,g)$ can be characterised as follows
\begin{equation}
  \lambda_1 = \inf_{0\neq u \in W_0^{1,2}(M) }\Bigg\{ \frac{\int_M | \nabla u |_g^2 \ d \mu_g}{\int_M | u|^2_g \ d\mu_g} \ \  \Big| \ \ u \neq 0, \ \ u \in W_0^{1, 2}(M,g) \Bigg\},
 \end{equation}
satisfying $\int_M u d \mu_g = 0 $, where $ W_0^{1, 2}(M,g)$ is the completion of $C^\infty_0(M,g)$ with respect to the sobolev norm
$$\|u\|_{W^{1,2}} = \Big(\int_M |u|^2 d\mu + \int_M | \nabla u|^2 d\mu \Big)^{\frac{1}{2}}.$$
Note that the $p$-Laplacian $\Delta_{p,g}, \ (1\leq p \leq \infty)$ is a natural generalisation of the Laplace-Beltrami operator.

   In this section, we consider the eigenvalues of the Laplace-Beltrami operator under the geometric flow, assuming the least eigenvalue $\lambda = \lambda(t)$ is a function of time only. Next we discuss the proof of Theorem \ref{thm311}.
   
   Let $M$ be a closed Riemannian manifold and $g(t)$ evolve by the generalized geometric flow (\ref{eq11}) in the interval $0\leq t \leq T.$ Let $\varphi(t) = \varphi(t, x)$ be the corresponding eigenfunction to the first nonzero eigenfunction $\lambda(t) = \lambda(t, \varphi)$ of $ \Delta_{2, g(t)}= \Delta$, i.e, 
  \begin{equation}\label{eq21}
  - \Delta \varphi(t, x) = \lambda(t) \varphi(t, x).
  \end{equation}
Taking derivative with respect to time, we have 
$$ - \Big(\frac{\partial}{\partial t} \Delta \Big) \varphi(t, x) - \Delta   \frac{\partial}{\partial t} \varphi (t, x) = \Big(\frac{d}{d t} \lambda(t) \Big)  \varphi(t, x)  + \lambda(t) \frac{\partial}{\partial t} \varphi(t, x). $$
Multiplying the above by $\varphi(t, x)$ and integrate with respect to the volume measure on $M$, we have 
$$ - \int_M \varphi \Big(\frac{\partial}{\partial t} \Delta \Big) \varphi \ d\mu - \int_M \varphi \Delta   \frac{\partial}{\partial t} \varphi \ d\mu  = \frac{d}{d t} \lambda(t) \int_M   \varphi^2 \ d\mu  + \lambda(t) \int_M \varphi  \frac{\partial}{\partial t} \varphi \ d\mu.$$  
Notice that  by the application of integration by parts and (\ref{eq21})
$$ - \int_M \varphi \Delta   \frac{\partial}{\partial t} \varphi \ d\mu  =  \lambda(t) \int_M \varphi  \frac{\partial}{\partial t} \varphi \ d\mu,$$
then, we arrive at 
 \begin{equation}\label{eq22}
   \frac{d}{d t} \lambda(t) \int_M   \varphi^2 \ d\mu = -\int_M \varphi \Big(\frac{\partial}{\partial t} \Delta \Big) \varphi \ d\mu.  
\end{equation}
Use the evolution of the Laplacian under the geometric flow, (i.e., (4) of Lemma  (\ref{lem211})), so that we have 
\begin{align*}
  \frac{d}{d t} \lambda(t) \int_M   \varphi^2 \ d\mu & = - 2  \int_M  h_{ij}  \nabla^i \nabla^j \varphi\   \varphi\ d\mu - 2 \int_M \langle div\ h , \nabla \varphi \rangle \varphi \ d\mu \\
  \displaystyle & \hspace{2cm} + \int_M \langle \nabla \mathcal{H}, \nabla \varphi \rangle \varphi \ d\mu.
  \end{align*}  
We express the first  and the last terms of the above as follows 
\begin{align*}
- 2  \int_M \varphi\  h_{ij}  \nabla^i \nabla^j \varphi\    d\mu & = \int_M \nabla^i ( 2 \varphi h_{ij} ) \nabla^j \varphi \ d\mu \\
\displaystyle & =  2 \int_M ( \nabla^i \varphi h_{ij} + \varphi \nabla_i h_{ij} ) \ \nabla^j \varphi \ d\mu \\
\displaystyle & =  2 \int_M h_{ij} \nabla^i \varphi \nabla^j \varphi \ d\mu + 2 \int_M \varphi \langle div\ h, \nabla \varphi \rangle \ d\mu
\end{align*}  
and   
\begin{align*}
\int_M  \langle  \nabla \mathcal{H}, \nabla \varphi \rangle \varphi \ d\mu & =  - \int _M \mathcal{H} \nabla_i(\varphi \nabla_j \varphi) \ d\mu\\
\displaystyle & =  -  \int_M \mathcal{H} \nabla^i \varphi \nabla^j \varphi \ d\mu -   \int_M \mathcal{H}  \varphi \Delta \varphi \ d\mu \\
\displaystyle & =   -  \int_M \mathcal{H} | \nabla \varphi |^2 \ d\mu + \lambda(t) \int_M  \mathcal{H} \varphi^2 \ d\mu.
\end{align*}  
Putting these together we have 
$$   \frac{d}{d t} \lambda(t) \int_M   \varphi^2 \ d\mu =  \lambda(t) \int_M  \mathcal{H} \varphi^2 \ d\mu -   \int_M \mathcal{H} | \nabla \varphi |^2 \ d\mu  + 2  \int_M h_{ij} \nabla^i \varphi \nabla^j \varphi \ d\mu. $$
Hence, we have proved Theorem \ref{thm311} on the evolution of $\lambda_1$.

In some applications $\mathcal{H}$ may be required to be a constant or  bounded by a constant. Thus when one takes $h_{ij}$ to be Ricci curvature tensor one is talking about a manifold with constant scalar curvature. In this situation we have the following as a corollary.
\begin{corollary}
Let $(M, g(t))$ be a closed manifold evolving by the geometric flow (\ref{eq11}). Let $\lambda(t)$ be the first nonzero eigenvalue of the Laplacian $\Delta_{g(t)}$. Then if $\mathcal{H} \geq C > 0$ in $M \times  [t_0, t ]$ for some uniform constant $C$ along the flow we have 
\begin{equation*}
   \frac{d}{d t} \log  \lambda(t)  \geq  C 
\end{equation*}
    and 
\begin{equation*}
  \lambda(t) \geq \lambda(t_0) e^{C(t - t_0)} \ for \ t > t_0.
  \end{equation*}
\end{corollary}

\section{Nonlinear eigenvalue problem for $p$-Laplacian}\label{sec4}
In this section, we consider the nonlinear eigenvalue problem  
\begin{equation}\label{eq41}
 \Delta_{p, g} u = - \lambda |u|^{p-2} u, \hspace{1cm}  u \neq 0 \ \ \ on \ M \times [0,T]
\end{equation}
with the normalization condition $\int_M |u|^p\ d\mu =1$.
We want to derive general evolution for the $p$-eigenvalues (eigenvalues of $\Delta_{p,g}$) and show that $\lambda_{p,1}$ is monotone  on metrics evolving by the geometric flow. In order to do these we need to calculate time evolution for $\lambda_{p,1}$ and its corresponding eigenfunction. But unfortunately, we do not know whether $\lambda_{p,1}$ or its corresponding eigenfunction $(p\neq2)$ is $C^1$-differentiable or not along the flow. So a similar approach to the one in \cite{[Ca07]} (see also Wu \cite{[Wu11]}, Wu, Wang and  Zheng \cite{[WWZ]} and Zhao \cite{[Zha13]}) will be used to avoid this difficulty. Precisely, let $(M, g_{ij}(t)), t \in [0, T]$ be a smooth compact Riemannian manifold evolving by the flow (\ref{eq11}). 
Define a genral smooth function as follows
 \begin{equation}
 \lambda_{p,1}(u(t), t) : = - \int_M u(t) \Delta_p u(t) d\mu_{g(t)} =  \int_M |\nabla u(t) |^pd\mu_{g(t)},
 \end{equation}
where $u(t)$ is a smooth  function satisfying  the normalisation  condition
 \begin{equation}\label{eqdef}
 \int_M|u(t)|^p d\mu_{g(t)} = 1, \ \ and \  \ \int_M |u(t)|^{p-2} u(t)d\mu_{g(t)} = 0.
\end{equation} 
  
By this we claim that there exists a smooth function $u(s)$ at time $t=s \in [0,T]$ satisfying (\ref{eqdef}). To see this claim, we first assume that at time $t=s$, $u(s)$ is the eigenfunction corresponding to $\lambda_{p,1}(s)$ of $\Delta_{p.g(s)}$ which implies
 \begin{equation*}
 \int_M|u(s)|^p d\mu_{g(s)} = 1, \ \ and \  \ \int_M |u(s)|^{p-2} u(s)d\mu_{g(s)} = 0.
\end{equation*}
Then we consider the following smooth function
\begin{equation}\label{eqdef1}
h(t) = u(s) \Bigg( \frac{|g(s)|}{|g(t)|}\Bigg)^{\frac{1}{2(p-2)}}
\end{equation} 
under the flow $g(t)$. We normalize this smooth function
\begin{equation}\label{eqdef2}
u(t) = \frac{h(t)}{\Big( \int_M |h(t)|^p d\mu_{g(t)} \Big)^{\frac{1}{p}}}
\end{equation} 
under the flow $g(t)$. 
By  (\ref{eqdef2}) we can easily check that $u(t)$ satisfies (\ref{eqdef}). 
 Note that in general $\lambda_{p,1}(u,t)$ is not equal to $\lambda_{p,1}(t)$. But at time $t=s$, if $u(s)$ is the eigenfunction of the first eigenvalue $\lambda_{p,1}(s)$, then we conclude that
$$\lambda_{p,1}(u(s),s) = \lambda_{p,1}(s).$$ 

Notice that the normalisation condition implies
 \begin{equation}
\frac{\partial} {\partial t}\Big(  \int_M |u|^p d \mu \Big)\Big|_{t=s} = 0,
\end{equation}     
which by direct computation (at $t=s$) yields the following
 \begin{align*}
\frac{\partial} {\partial t}\Big(  \int_M |u|^p d \mu \Big) = 0 &= \frac{\partial}
 {\partial t}\Big(  \int_M |u|^{p-1} u d \mu \Big)\\
\displaystyle& = \int_M(p-1)  |u|^{p-2} u \frac{\partial u} {\partial t} d\mu + 
 \int_M |u|^{p-1} \frac{\partial } {\partial t} (u d\mu).
\end{align*}   
By this it holds that 
\begin{equation}\label{eq47}
\int_M |u|^{p-2} u \Big((p-1) \frac{\partial u} {\partial t} d\mu +  \frac{\partial }
 {\partial t} ( u d\mu) \Big) = 0.
\end{equation}

We now present a proposition

\begin{proposition}
Let $g(t)$ be a smooth solution of the flow (\ref{eq11}) on $M$. Let $\lambda_{p,1}(t)$ be the
  first eigenvalue of the $p$-Laplacian under (\ref{eq11}). Assume $u(s)$ is the corresponding eigenfunction of $\lambda_{p,1}(t)$ at time $t=s \in [0, T]$, that is,
\begin{equation*}
\Delta_{p,g(s),p} u(s) = - \lambda_{p,1}(s) |u(s)|^{p-2} u(s). 
\end{equation*}  
Let $\lambda_{p,1}(u,t)$ be a smooth function defined by 
 \begin{equation}
 \lambda_{p,1}(u(t), t) : = - \int_M u(t) \Delta_p u(t) d\mu_{g(t)}.
 \end{equation}  
Then we have   
\begin{equation}\label{eq5.9}
\left. \begin{array}{l}
\displaystyle\frac{d}{dt} \lambda_{p,1}(u,t)\Big|_{t=s} = \lambda_{p, 1}(s) \int_M S_g | u(s)|^p
  d\mu_{g(s)} -  \int_M S_g | \nabla u(s)|^p d\mu_{g(s)} \\ 
\displaystyle \hspace{3cm}+ p  \int_M  | \nabla u(s)|^{p-2} \mathcal{S}_{ij} u_i u_j   d\mu_{g(s)},
\end{array} \right. 
\end{equation}
where $u(t)$ is any smooth function satisfying (\ref{eqdef}) such that at any time $ t = s \in[0, T]$, $u(s)$ is the eigenfunctin for $\lambda_{p,1}(s)$. Here $u_i$ denotes covariant derivative of $u$ with respect to the Levi-Civita connection of $g(t)$ at time $t=s$.
\end{proposition}
    
Since the evolution formula (\ref{eq5.9}) does not depend on the time derivative of the eigenfunction $u$, we have that
$$\frac{d}{dt} \lambda_{p,1}(u(t), t)=\frac{d}{dt}\lambda_{p,1}(t)$$
at some time $t=s$.
 
We are now set to prove a theorem about the evolution,  monotonicity and differentiability (Theorem \ref{thm411}) of the first eigenvalue of the $p$-Laplacian under the geometric flow.  
Clearly, we can now set   
 \begin{equation}
 \lambda_{p,1}(t) = \lambda_{p,1}(u(t), t)  = -  \int_M u(t,x) \Delta_p u(t,x) d\mu_{g(t)}.
\end{equation}
The evolution of $\lambda_{p,1}$ then follows
\begin{equation}
\frac{d}{dt}\lambda_{p,1}(t) = \frac{d}{dt} \lambda_{p,1}(u(t), t)  = - \frac{\partial} {\partial t} \int_M u(t,x)
 \Delta_p u(t,x) d\mu_{g(t)},
\end{equation}

\subsection*{Proof of Theorem \ref{thm411}}
\begin{proof}
The proof also follows by direct computation using evolution of quantities in Lemmas \ref{lem211} and \ref{lem212}.
Denote $Z:= |\nabla u|^{p-2}$, then working in local orthonormal coordinates we have 
 \begin{align*}
\frac{\partial} {\partial t}   \int_M u \Delta_p u d\mu &=  \frac{\partial} {\partial t}   \int_M  g^{ij} \nabla_i[Z \nabla_j u ] u d\mu  \\
\displaystyle& =   \frac{\partial} {\partial t} \int_M   \Big( g^{ij} \nabla_i Z \nabla_j u + Z \Delta u \Big) u d\mu \\
\displaystyle& =  \int_M \frac{\partial} {\partial t}   \Big( g^{ij} \nabla_i Z \nabla_j u + Z \Delta u \Big) u d\mu + \int_M \Delta_p u \frac{\partial} {\partial t}   (u d\mu)\\
\displaystyle& =: I + II.
\end{align*}  
By the evolution of $\Delta_{p,g}$ in Lemma \ref{lem212} we have 
\begin{align*}
 I & =  \int_M  2 h^{ij} \nabla_i ( Z \nabla_j u ) + g^{ij} \nabla_i ( Z_t \nabla_j u ) + g^{ij} \nabla_i (Z \nabla_j u_t ) \\
   \displaystyle &\hspace{1cm} +   \int_M Z \Big\{2 \langle \ div \ h , \nabla u \rangle -   \langle \nabla \mathcal{H}  , \nabla u \rangle \Big\} u d\mu.
 \end{align*}
Using integration by parts on the second and third terms of the last integral we have 
\begin{align*}
 I & =  \int_M  2 h^{ij} \nabla_i ( Z \nabla_j u ) -   \int_M g^{ij}  Z_t \nabla_i u \nabla_j u  - \int_M   Z g^{ij}  \nabla_i u \nabla_j u_t ) \\
   \displaystyle &\hspace{1cm} + \int_M   Z \Big\{2 \langle \ div \ h , \nabla u \rangle -   \langle \nabla \mathcal{H}  , \nabla u \rangle \Big\} u d\mu.
 \end{align*}
Therefore we have after using the evolution $Z_t$ from Lemma \ref{lem212}
\begin{equation}\label{eq411}
\left. \begin{array}{l}
\displaystyle\frac{\partial} {\partial t}   \int_M u \Delta_p u d\mu = \int_M  2 h^{ij} \nabla_i ( Z \nabla_j u ) ud\mu -(p-2) \int_M |\nabla u|^{p-2} h^{ij} \nabla_i u \nabla_j u d\mu\\
\displaystyle \hspace{3cm} -(p-1) \int_M |\nabla u|^{p-1} g^{ij} \nabla_i u \nabla_j u_t d\mu \\
 \displaystyle \hspace{3cm} +   \int_M Z \Big\{2 \langle \ div \ h , \nabla u \rangle -   \langle \nabla \mathcal{H}  , \nabla u \rangle  \Big\} u d\mu \\
\displaystyle \hspace{3cm}  +  \int_M \Delta_p u \frac{\partial} {\partial t}   (u d\mu).
\end{array} \right. 
\end{equation}
Computing the first and the third terms on the right hand side (RHS for short) of (\ref{eq411})as follows
\begin{align*}
\int_M 2 h^{ij} \nabla_i ( Z \nabla_j u ) u d\mu  &= - 2 \int_M  \nabla_i ( h^{ij} u) Z \nabla_j u u d\mu \\
\displaystyle&=- 2 \int_M |\nabla u|^{p-2} h^{ij} \nabla_i u \nabla_j u d\mu - 2 \int_M Z \langle div h, \nabla u \rangle u d\mu.
\end{align*}
\begin{align*}
 -(p-1) \int_M |\nabla u|^{p-1} g^{ij} \nabla_i u \nabla_j u_t d\mu &= (p-1) \int_M g^{ij}\nabla_i( |\nabla u|^{p-1} \nabla_i u ) u_t \\
 \displaystyle&=  (p-1) \int_M \Delta_p u u_t d\mu
 \end{align*}
Putting these back into (\ref{eq411}) we have
\begin{align*}
\frac{\partial} {\partial t}   \int_M u \Delta_p u d\mu  &= -p \int_M |\nabla u|^{p-2} h^{ij} \nabla_i u \nabla_j u d\mu - \int_M Z \langle \nabla \mathcal{H}  , \nabla u \rangle  u d\mu \\
 \displaystyle& \ \ \ +  (p-1) \int_M \Delta_p u \ u_t d\mu +  \int_M \Delta_p u \frac{\partial} {\partial t}   (u d\mu).
  \end{align*}
  Using the integrability condition (\ref{eq47}) and the definition of $\Delta_p u$ in (\ref{eq41}), the last two terms on the RHS of the above equation vanish and we then arrive at 
 \begin{align}\label{eq412}
\frac{\partial} {\partial t}   \int_M u \Delta_p u d\mu  = -p \int_M |\nabla u|^{p-2} h^{ij} \nabla_i u \nabla_j u d\mu - \int_M Z \langle \nabla \mathcal{H}  , \nabla u \rangle  u d\mu.
 \end{align}
The next is to compute the second term on the RHS of the last equality using integration by parts as follows
\begin{align*}
- \int_M Z \langle \nabla \mathcal{H}  , \nabla u \rangle  u d\mu & = \int_M \mathcal{H} \nabla_i(Z \nabla_j u \ u ) d\mu  \\
\displaystyle& = \int_M \mathcal{H} \nabla_i(Z \nabla_j u ) u d\mu + \int_M \mathcal{H} Z \nabla_i u \nabla_j u d\mu \\
\displaystyle& = \int_M \mathcal{H} \Delta_p u\ u d\mu + \int_M \mathcal{H} |\nabla u|^{p-2} | \nabla u|^2 d\mu \\
\displaystyle& = -  \lambda_{p, 1}(t) \int_M \mathcal{H} |u|^p d\mu + \int_M \mathcal{H} |\nabla u|^p d\mu. 
 \end{align*}
Putting this into (\ref{eq412}) we obtain (\ref{eq48}) at once. Hence, the first part of the theorem is proved. Using the condition (\ref{eq49}) in (\ref{eq48}) we have the monotonicity formula (\ref{eq49}) with the condition $\mathcal{H} \geq 0$.
\end{proof}
Clearly, when $p=2$ we have $\Delta_2 = \Delta_g$, the usual Laplace-Beltrami operator. Also $\lambda_{p, 1} = \lambda_1$, the first eigenvalue of $\Delta_g$ and the corresponding eigenfunction are smoothly differentiable. Then Theorem \ref{thm311} reduces to a corollary. This further explains that the $p$-Laplacian is a nonlinear generalisation of Laplace-Beltrami operator.

Integrating both sides of (\ref{eq410}) from $t_1$ to $t_2$ on a sufficiently small time interval $ t_1 \leq t \leq t_2$, $ t \in [0, T]$, we then obtain \ref{eq514}. 

\begin{proof} of Corollary \ref{cor4}.
By the definition of $ \lambda_{p, 1}$ in (\ref{eq13}) and the normalization condition $\int_M |u|^p d\mu =1$ we know that
$$ \lambda_{p, 1} = \int_M  |\nabla u|^p d\mu.$$
Then  (\ref{eq410}) reduces to 
\begin{equation}\label{eq415}
\frac{d}{dt} \lambda_{p, 1}(t) \geq  \alpha p \psi(t)  \lambda_{p, 1}(t)
\end{equation}
from which (\ref{eq414}) follows by integrating on the interval $[t_1, t_2]$ with $ t_1, t_2 \in [0, T].$

\end{proof}

Note that both $\lambda_1(t)$ and $\psi(t)$ are functions of time only.
Setting
$$\psi(0) = \mathcal{H}_{min}(0) = \psi_0,$$
we can evaluate
\begin{align*}
\int_{t_1}^{t_2} \psi(t) dt &= \int_{t_1}^{t_2} \Big(  \frac{1}{\psi^{-1}_0 - \frac{2}{n}t} \Big) dt \\
\displaystyle & = - \frac{n}{2} \log ( \psi^{-1}_0 - \frac{2}{n}t) \Big|_{t_1}^{t_2} = \log \Bigg( \frac{ \psi^{-1}_0 - \frac{2}{n}t_1}{ \psi^{-1}_0 - \frac{2}{n}t_2} \Bigg)^{\frac{n}{2}}.
\end{align*}
Therefore integrating both sides of (\ref{eq415}) from $t_1$ to $t_2$ and together with the above equality, yields 
\begin{equation}
\log  \frac{\lambda_{p,1}(t_2)}{ \lambda_{p,1}(t_1)} \geq \log \Bigg( \frac{ \psi^{-1}_0 - \frac{2}{n}t_1}{ \psi^{-1}_0 - \frac{2}{n}t_2} \Bigg)^{\frac{\alpha np}{2}}
\end{equation}
for any time $t_1 < t_2$ sufficiently close to $t_2$. By this we have 
\begin{equation*}
\lambda_{p,1}(t_2) \cdot \Big( \psi^{-1}_0 - \frac{2}{n}t_2 \Big)^{\frac{\alpha np}{2}} \geq \lambda_{p,1}(t_1) \cdot \Big(\psi^{-1}_0 - \frac{2}{n}t_1\Big)^{\frac{\alpha np}{2}}.
\end{equation*}
Then $\lambda_{p,1}(t) \cdot ( \psi^{-1}_0 - \frac{2}{n}t)^{\frac{\alpha np}{2}}$ is nondecreasing along the geometric flow. Hence we conclude this section with the following.

\begin{theorem}\label{thm44}
With the assumption of Theorem \ref{thm411}.  The following quantity 
\begin{equation*}
\lambda_{p,1}(t) \cdot \Big( \psi^{-1}_0 - \frac{2}{n}t \Big)^{\frac{\alpha np}{2}}
\end{equation*}
is nondecreasing  and $\lambda_{p,1}(t) $ is differentiable almost everywhere along the geometric flow (\ref{eq11}).
\end{theorem}
The last theorem has been proved using a different method under the Ricci flow  in \cite{[WWZ]}.

\subsection*{The differentiability of $\lambda_{p,1}(t)$}
 Since $\lambda_{p,1}(t)$ is nondecreasing on the time interval $[0,T)$ (under  curvature assuption of the theorem), by the classical Lebesgue's theorem [A. Mukherjea and K. Pothoven: Real and Functional Analyisis (Chap 4)], it is easy to see that $\lambda_{p,1}(t)$ is differentiable almost everywhere, 

\begin{remark}
Our proofs of the first eigenvalue evolution and monotonicity do not use any differentiability of the first eigenvalue $\lambda_{p,1}(t)$ or its corresponding eigenfunction $u(t,x)$ of the $p$-Laplacian under $(RH)_\alpha$-flow. In fact, it is not known whether they are differentiable in advance.

It would be interesting to find out whether the corresponding first eigenfunction of the $p$-Laplacian is a $C^1$-differentiable function with respect to $t$-variable along $(RH)_\alpha$-flow.
\end{remark}


\section{Examples of geometric flows} 
In this section, we give some examples of geometric flows where our results are valid. We remark that in these cases the error term $\mathcal{D}(X)$ and the quantity $\beta - \Delta \mathcal{H}$ are nonnegative. More examples can be found in \cite[Section 2]{[Mu10]}.

\subsection{Hamilton's Ricci  flow \cite{[Ha82]}}
Let $(M, g(t))$ be a solution to the Hamilton's Ricci flow
\begin{equation}
\partial_t g_{ij}(t, x) = - 2 R_{ij}.
\end{equation}
This is the case where $h_{ij} =R_{ij}$ is the Ricci tensor and $\mathcal{H} =R$ is the scalar curvature on $M$. Here, the scalar curvature evolves by 
$$\partial_t R = \Delta R + 2 |R_{ij}|^2.$$
 By twice contracted second Bianchi identity $g^{ij} \nabla_i R_{jk} = \frac{1}{2} \nabla_k R$, which implies 
 $$  2 \langle \ div \ h , \nabla f \rangle -   \langle \nabla \mathcal{H}   , \nabla f \rangle = 0, $$
  the quantity $\mathcal{D}(X)$ vanishes identically and $\beta -\Delta R  \equiv 0$. Note that the positivity of curvature is preserved along the Ricci flow \cite{[Ha82]}. The evolution equation and monotonicity formula  for the first eigenvalue follow easily (see \cite{[Ca07]}-\cite{[GPT]} and \cite{[Li07a]}-\cite{[Ma06]}). A fundamental result here is Perelman's paper \cite{[Pe02]}, where he defines his  $\mathcal{F}$-energy 
  \begin{equation}
 \mathcal{F}(g_{ij}(t), u(t)) = \int_{M} (4 |\nabla u|^2  + R u^2) d \mu \ \ \ \ \ \  with \ \ \  \int_M  u^2 d \mu = 1.  
 \end{equation}
 and proves that it is monotonically nondecreasing. He also defines
  \begin{equation}
 \lambda_1 ( g_{ij}) = \inf \Big\{ \mathcal{F} (g_{ij}, f ) : f \in C_c^\infty (M) , \int_M  e^{ - f} d \mu = 1 \Big\},
 \end{equation}
with $ \lambda_1 ( g_{ij})$ (being the least eigenvalue of the geometric operator $- 4\Delta + R$) and its corresponding eigenfunction $u=e^{ - f} $ satisfying the eigenvalue problem
$$ - 4 \Delta u + R u =  \lambda_1 ( g_{ij}) u.$$
He shows that monotonicity of $ \lambda_1 ( g_{ij})$ follows from that of  $\mathcal{F}$.

\subsection{Ricci-harmonic map flow \cite{[Mu12]} }
Let $(M, g)$ and $(N, \xi)$ be compact (without boundary) Riemannian manifolds of dimensions $m$ and $n$ respectively.
Let a smooth map $ \varphi: M  \rightarrow N$ be a critical point of the Dirichlet energy integral $E( \varphi) = \int_M | \nabla  \varphi|^2 d \mu_g$, where $N$ is isometrically embedded in $\mathbb{R}^d, \ d \geq n,$ by the Nash embedding theorem. The configuration $(g(x ,t),  \varphi(x, t)), t \in [0, T)$ of  a one parameter family of Riemannian metrics $g(x, t)$ and a family of smooth maps $ \varphi(x, t)$ is defined to be Ricci-harmonic map flow if it satisfies the coupled system of nonlinear parabolic equations 
\begin{equation}\label{eq514}
\left \{ \begin{array}{l}
\displaystyle \frac{\partial}{\partial t} g(x ,t) = - 2 Rc(x ,t) + 2 \alpha \nabla  \varphi(x ,t) \otimes \nabla  \varphi(x ,t)
\\ \ \\
\displaystyle \frac{\partial}{\partial t}  \varphi(x ,t) = \tau_g  \varphi(x ,t),
\end{array} \right.
\end{equation}
where $ Rc(x ,t)$ is the Ricci curvature tensor for the metric $g$, $\alpha(t) \equiv \alpha > 0$ is  a time-dependent coupling constant, $\tau_g  \varphi$ is the intrinsic Laplacian of $ \varphi$, which denotes the tension field of map $ \varphi$ and $\nabla \varphi \otimes \nabla  \varphi = \varphi^* \xi$ is the pullback of the metric $\xi$ on $N$ via the map $\varphi$. See List \cite{[Li08]} when the target manifold is one dimensional.
Here $h_{ij} = R_{ij} - \alpha \partial_i \varphi \partial_j \varphi =:S_{ij}$, \  $\mathcal{H}= R - \alpha |\nabla \varphi|^2 =: S$ and 
\begin{equation}
 \partial_t S = \Delta S + 2|S_{ij}|^2 + 2 \alpha |\tau_g \varphi|^2 - 2 \dot{\alpha} |\nabla \varphi|^2.
 \end{equation}
Using the twice contracted second Bianchi identity, we have 
\begin{equation}
(g^{ij} \nabla_i S_{jk} - \frac{1}{2} \nabla_k S) X_j = - \alpha \tau_g \varphi \nabla_j\varphi X_j.
\end{equation} 
A straightforward computation gives
\begin{equation}\label{eq44}
2 (Rc - \mathcal{S}) (\nabla u, \nabla u) = 2 \alpha  | \nabla \varphi|^2 | \nabla u |^2 .
\end{equation} 
Then, $ \mathcal{D}(X) = 2 \alpha | \tau_g \varphi - \nabla_X \varphi|^2 - 2 \dot{\alpha}  |\nabla \varphi|^2 $ and $\beta - \Delta S =  2 \alpha | \tau_g \varphi|^2 - \dot{\alpha}  |\nabla \varphi|^2 $ for all $X$ on $M$. Thus both $\mathcal{D}$ and $\beta - \Delta S$ are nonnegative as long as $\alpha(t)$ is nonincreasing in time. The first author has considered this in \cite{[Ab15]}. See also \cite{[Li10]} for the monotonicity of the first eigenvalue of Laplace-Betrami operator and versions of Perelman's entropy under the Ricci-harmonic map flow.

\subsection{Lorentzian mean curvature flow}

Let $M^n(t) \subset L^{n+1}$ be a family of space-like hypersurfaces in ambient 
 Lorentzian manifold evolving by Lorentzian mean curvature flow
 $$\partial_t F(t, \cdot) = \Pi(t, \cdot) \nu(t, \cdot)$$
for $(t, \cdot) \in [0,T] \times M$, where $F(t, \cdot)$ is the position of $M^n$ in $L^{n+1}$ satisfying $F(0,\cdot)=F_0(\cdot)$. Here $\nu(t, \cdot)$ and $\Pi(t, \cdot)$ 
 are respectively the outer normal vector and mean curvature at the point $F(t, \cdot).$ Then, the induced metric evolves by 
$$\partial_tg_{ij} = 2 H \Pi_{ij}, $$
where $\Pi_{ij}$ denotes the  components of the second fundamental form $\Pi$ on $M$ and $H = g^{ij} \Pi_{ij}$ denotes the mean curvature of $M$. In this case 
$h_{ij} = - H \Pi_{ij}$ and $\mathcal{H} = - H^2$. Letting $ \widetilde{Rc}$ and $\widetilde{Rm}$ denote the Ricci and Riemman curvature tensor of $L^{n+1}$ respectively, we have the Gauss equation
$$ R_{ij} = \widetilde{R}_{ij} - H \Pi_{ij} + \Pi_{il} \Pi_{lj} + \widetilde{R}_{i0j0},$$
the Codazzzi equation
$$\nabla_i \Pi_{jk} - \nabla_k \Pi_{ij} = \widetilde{R}_{0jki},$$
the evolution equation
$$\partial_tH = \Delta H - H (|\Pi|^2 + \widetilde{Rc}(\nu, \nu)$$ 
and
$$\beta - \Delta H = 2 H^2 |\Pi|^2 + |\nabla H|^2 + 2 H \widetilde{Rc}(\nu, \nu).$$
See the explicit forms  of the Gauss  and  the Codazzi equations 
 for $L^{n+1} = \mathbb{R}^{n+1}$ in \cite{[Hu84]}. Combining the above equation we obtain the quantity
\begin{equation}
 \mathcal{D}(X) = 2|\nabla H - \Pi(X, \cdot)|^2 + 2 \widetilde{Rc}(H \nu - X, 	 H \nu - X) + 2 \langle \widetilde{Rm}(X, \nu)\nu, X \rangle,
 \end{equation}
 where $\nu$ denotes future-oriented timelike unit normal vector on $M$. Obviously both $\mathcal{D}(X)$ and $\beta -\Delta H$ are nonnegative when assuming nonnegativity on sectional curvature of $L^{n+1}$. See \cite{[Zha13]} for the evolution and monotonicity of the  first eigenvalue of $p$-Laplace operator under the  $m^{th}$ powers of the mean curvature flow. See also Huisken \cite{[Hu84]}, the second author's paper \cite{[Mao1]} and \cite{[GLW16],[Zha12]} for related results.

\subsection{The Yamabe flow}
This is the case when $h_{ij} = \frac{1}{2}Rg_{ij}$, where $R$ is the scalar curvature of the metric. Yamabe flow is then the following evolution equation
\begin{equation}\label{eq518}
\left. \begin{array}{l}
\displaystyle \frac{\partial}{\partial t} g_{ij}(x ,t) = - R(x, t) g_{ij}(x ,t),  \ \ (x, t) \in M \times [0, T]
\\ \ \\
g_{ij}(x ,0)  =  g_0(x)
\end{array} \right.
\end{equation}
as introduced by R. Hamilton who first establishes the existence of its unique solution for all time and shows that the metric $g(t)$ approaches constant as $t \rightarrow \infty$. His proof is done for volume preserving flow
\begin{equation}
\frac{\partial}{\partial t} g_{ij}(x ,t) = (r(t)- R(x, t)) g_{ij}(x ,t),  \ \ (x, t) \in M \times [0, \infty),
\end{equation}
with $r(t) = Vol^{-1}(g(t))\int_M R d\mu$ is the average of scalar curvature for the metric in a conformal class. For more on the global existence and convergence of (\ref{eq518}) see Chow \cite{[Ch92]} and Ye \cite{[Ye94]}. Note that under the Yamabe flow the volume measure evolves as $\partial_t d\mu =  n/2R d\mu$ and the normalization condition, $\partial_t(\int_M |u|^p d\mu) = 0$, implies 
$$ \int_M p |u|^p u u_t d\mu - \frac{n}{2} \int_M R u d\mu = 0.$$
Here, the evolution of scalar curvature is given \cite{[Ch92]} as 
\begin{equation}
\frac{\partial}{\partial t} R = (n-1) \Delta R + R^2
\end{equation}
and by the strong maximum principle
$$R(x, t) \geq \psi(t) = \frac{R_{min}(0)}{1 - R_{min}(0)  t} $$
for all t.
We can also compute 
\begin{align*}
\langle 2 div\ h - \nabla \mathcal{H}, \nabla f \rangle &= \langle 2 g^{ij} \nabla_i( \frac{1}{2} Rg_{ij}) - \frac{n}{2} \nabla_k R, \nabla f \rangle \\
\displaystyle & = \frac{2n-n}{2} \langle \nabla R, \nabla f \rangle
\end{align*}
and 
$$\beta - \Delta \mathcal{H} = \frac{n(n-2)}{2} \Delta R.$$
These imply that the quantitty $\mathcal{D}(\nabla f)$ is nonnegative on the Einstein tensor
$$E_{ij} = R_{ij} - \frac{1}{2} R_{ij} \geq 0.$$
Hence our results hold. See \cite[Section 7]{[WWZ]}, where with assumption  that $p \geq n$ and $R \geq 0$, they prove that $\lambda_{p, 1}$ is strictly increasing and differentiable almost everywhere along the Yamabe flow. Yamabe flow coincides with the Ricci flow on Riemann surfaces.

\subsection*{Acknowledgements}
The first author was partially supported by Nigeria Tetfund research grant through OSCOTECH, Esa-Oke. The second author was partially supported by the NSF of China (Grant No. 11401131). The first version of the paper was completed during the first author's research visit to  African Institute for Mathematical Sciences (AIMS), Senegal in February 2016. He therefore thanks AIMS-Senegal for finacial suport during the visit. He also wishes to thank the Research chair, Prof. Moustapha Fall for the hospitality and useful discussions which helped in improving this paper.


\end{document}